\documentclass{amsart}
\usepackage[all]{xy}

\newtheorem{theorem}{Theorem}[section]
\newtheorem{lemma}[theorem]{Lemma}

\theoremstyle{definition}
\newtheorem{definition}[theorem]{Definition}
\newtheorem{example}[theorem]{Example}
\newtheorem{remark}[theorem]{Remark}

\renewcommand{\P}{{\mathbb{P}}}
\newcommand{\Z}{{\mathbb{Z}}}
\newcommand{\Fc}{{\mathcal{F}}}
\newcommand{\Oc}{{\mathcal{O}}}

\newcommand{\ed}{\widehat{E}}

\newcommand\res[1]{{\lower1pt\hbox{$|$}}_{\raise.5pt\hbox{${\scriptstyle #1}$}}}

\numberwithin{equation}{section}

\begin{document}

\title{Bezoutians and Tate Resolutions} 

\author{David A.\ Cox}
\address{Department of Mathematics and Computer Science, Amherst 
College, Amherst, MA 01002-5000}
\email{dac@cs.amherst.edu}


\begin{abstract}
This paper gives an explicit construction of the Tate resolution of
sheaves arising from the $d$-fold Veronese embedding of $\P^n$.  Our
description involves the Bezoutian of $n+1$ homogenous forms of degree
$d$ in $n+1$ variables.  We give applications to duality theorems,
including Koszul duality.
\end{abstract}

\date{\today}

\maketitle

\section{Introduction}

Given a finite dimensional vector space $W$ over a field $k$ with dual
$V$, a coherent sheaf $\Fc$ on $\P(W)$ gives a \emph{Tate resolution}
$T^\bullet(\Fc)$, which is a minimal bi-infinite exact sequence of
free graded $E =\bigwedge V$-modules
\[
\cdots \longrightarrow T^{-2}(\Fc) \longrightarrow
T^{-1}(\Fc) \longrightarrow T^{0}(\Fc) \longrightarrow
T^{1}(\Fc) \longrightarrow T^{2}(\Fc) \longrightarrow
\cdots. 
\]
These resolutions were introduced by Gel$\!$\'{}$\!$fand
\cite{gelfand} in 1984 and are part of the BGG correspondence
\cite{bernstein} from 1978.

The paper \cite{efs} gives an explicit formula for
$T^\bullet(\Fc)$, namely
\begin{equation}
\label{tp}
T^p(\Fc) = {\textstyle \bigoplus_i} \ed(i-p)\otimes_k H^i(\P(W),\Fc(p-i)),
\end{equation}
where $\ed = \mathrm{Hom}_k(E,k)= \bigwedge W$ as an $E$-module.  Also
note that $\deg(W) = 1$ since $\deg(V) = -1$ and that $\ed \simeq
E(-\dim(W))$ (noncanonically).

The maps $T^p(\Fc) \to T^{p+1}(\Fc)$ are less well understood.  For
the $i$th summand of $T^p(\Fc)$, the map to $T^{p+1}(\Fc)$ looks like
\[
\xymatrix@R=1pt{ 
\ed(i-p)\otimes_k H^i(\Fc(p-i)) \ar[r] \ar[ddr]
\ar[ddddddr]^(.3){\raise3pt\hbox{$\vdots$}}& 
\ed(i-p-1)\otimes_k H^i(\Fc(p+1-i)) \\ 
& \bigoplus \\ 
& \ed(i-p-2)\otimes_k H^{i-1}(\Fc(p+2-i)) \\ 
& \bigoplus \\ & \raise7pt\hbox{$\vdots$} \\ 
& \bigoplus \\ & \ed(-p)\otimes_k H^0(\Fc(p)),}
\]
where for simplicity we have omitted ``$\P(W)$'' in the cohomology
groups.  The horizontal map in this diagram is known from \cite{efs},
while the diagonal maps are more mysterious.  Examples of these
diagonal maps can be found \cite{efs,es}, and explicit descriptions of
certain diagonal maps in the toric context were given by Khetan in his
work \cite{khetan1,khetan2} on sparse determinantal formulas in
dimensions 2 and 3.

In this paper, we will use Bezoutians to describe the diagonal maps in
the Tate resolution for a particular choice of $\Fc$.  Let $S =
k[x_0,\dots,x_n]$ have the standard grading and let $W = S_d$ be the
graded piece in degree $d \ge 1$.  Thus $\dim(W) = \binom{n+d}{d}$.
Given any $\ell \in \Z$, the $d$-fold Veronese embedding
\[
\nu_d : \P^n \longrightarrow \P(W)
\]
gives the coherent sheaf
\[
\Fc = \nu_{d*}\Oc_{\P^n}(\ell)
\]
on $\P(W)$.  We will give an explicit construction of the Tate
resolution $T^\bullet(\Fc)$.

Since $\Oc_{\P(W)}(1)\res{\nu_d(\P^n)} = \nu_{d*}\Oc_{\P^n}(d)$, we
have
\[
H^i(\P(W),\Fc(j)) = H^i(\P^n,\Oc_{\P^n}(\ell + jd)).
\]
This cohomology group will be denoted $H^i(\ell + jd)$.  Using Serre
duality and standard vanishing theorems for line bundles on $\P^n$, we
also have
\[
H^i(\ell + jd) = \begin{cases} S_{\ell + jd} & i =
  0\\
S_{-n-1-(\ell+jd)}^* & i = n\\
0 & \text{otherwise,}
\end{cases}
\]
where $S_m$ is the graded piece of $S = k[x_0,\dots,x_n]$ in degree
$m$.

In the Tate resolution, it follows that
\begin{align*}
T^p(\Fc) &= \ed(-p)\otimes_k H^0(\ell+pd)\ {\textstyle\bigoplus}\
\ed(n-p)\otimes_k H^n(\ell+(p-n)d)\\ 
&= \ed(-p)\otimes_k S_{\ell+pd}\ {\textstyle\bigoplus}\ \ed(n-p)\otimes_k
S_{-n-1-(\ell+(p-n)d)}^*.
\end{align*}
To simplify the subscripts, we set $a = \ell + (p+1)d$ and $\rho =
(n+1)(d-1)$.  Then the description of $T^p(\Fc)$ becomes
\[
T^p(\Fc) = \ed(-p)\otimes_k S_{a-d}\ {\textstyle\bigoplus}\
\ed(n-p)\otimes_k S_{\rho - a}^*,
\]
and the map $T^p(\Fc) \to T^{p+1}(\Fc)$ has the following form:
\begin{equation}
\label{tpf}
\xymatrix@R=5pt{
\widehat{E}(n-p)\otimes_k S_{\rho - a}^* \ar[r]^(.45){\alpha_p}
\ar[ddr]^{\delta_p} & \widehat{E}(n-p-1)\otimes_k S_{\rho - a - d}^*\\
\bigoplus & \bigoplus\\
\widehat{E}(-p)\otimes_k S_{a-d} \ar[r]^{\beta_p} &
\widehat{E}(-p-1)\otimes_k S_{a}.}
\end{equation}
By \cite{efs}, the map 
\[
\beta_p \in \mathrm{Hom}_E(\ed(-p)\otimes_k
S_{a-d},\ed(-p-1)\otimes_k S_a)_0 \simeq
\mathrm{Hom}_k(W\otimes_k S_{a-d}, S_a)
\]
(the subscript ``0'' means graded $E$-module homomorphisms of degree
0) corresponds to multiplication $W\otimes_k S_{a-d} = S_d \otimes_k
S_{a-d} \to S_a$, and $\alpha_p$ similarly corresponds to the
natural map $W\otimes_k S^*_{\rho-a} \to S^*_{\rho-a-d}$ induced by
multiplication.

The diagonal map $\delta_p$ in \eqref{tpf} lies in
\begin{equation}
\label{hom}
\mathrm{Hom}_E(\ed(n-p)\otimes_k
S_{\rho-a}^*,\ed(-p-1)\otimes_k S_a)_0 \simeq
\mathrm{Hom}_k({\textstyle\bigwedge^{n+1}} W, S_{\rho-a}\otimes_k S_a).
\end{equation}
The map $\delta_p$ is not unique; hence our main result
(Theorem~\ref{main} below) will give one possible choice for the this
map.

We next recall the definition of the Bezoutian.  

\begin{definition}
\label{bezoutdef}
Consider the polynomial ring $k[x_0,\dots,x_n,y_0,\dots,y_n]$.    
\begin{enumerate}
\item For $f \in k[x_0,\dots,x_n]$ and $0 \le j \le n$, define
$\Delta_{j}(f)$ to be the polynomial
\[
\frac{f(y_0,\dots,y_{j-1},x_j,x_{j+1},\dots,x_n) -
  f(y_0,\dots,y_{j-1},y_j,x_{j+1},\dots,x_n)}{x_j-y_j}.
\]
\item The \emph{Bezoutian} of homogeneous polynomials $f_0,\dots,f_n
\in k[x_0,\dots,x_n]$ of degree $d$ is the $(n+1)\times(n+1)$
determinant
\[
\Delta = \det \Delta_{j}(f_i).
\]
\end{enumerate}
\end{definition}

\begin{remark} 
\label{bezoutprops}
Here are some observations about the Bezoutian of $f_0,\dots,f_n$.
\begin{enumerate}
\item Each $\Delta_{j}(f_i)$ is homogeneous of degree $d - 1$ in
$x_0,\dots,x_n,y_0,\dots,y_n$, so the Bezoutian is
homogeneous of degree $\rho = (n+1)(d-1)$ in these variables.
\item Writing $\Delta$ as a polynomial in the $y_i$'s with
  coefficients in $k[x_0,\dots,x_n]$, we obtain
\[
\Delta = \sum_{|\alpha|\le \rho} \Delta_\alpha(x) y^\alpha,
\]
where $\Delta_\alpha(x) \in S = k[x_0,\dots,x_n]$ has degree
$\rho-|\alpha|$. 
\item Under the natural bigrading of $k[x_0,\dots,x_n,y_0,\dots,x_n]$,
  the graded piece of $\Delta$ of bidegree $(\rho-a,a)$ is
\[
\Delta_{\rho-a,a} = \sum_{|\alpha| = a}\Delta_\alpha(x)y^\alpha.
\]
\item Recall the isomorphism $k[x_0,\dots,x_n,y_0,\dots,x_n] \simeq
S\otimes_k S$ given by $x_i \mapsto x_i\otimes1, y_i \mapsto 1\otimes
x_i$.  Since  $\Delta$ is multilinear and alternating in
$f_0,\dots,f_n$, the Bezoutian construction gives a linear map
\[
{\textstyle\bigwedge^{n+1}} S_d \longrightarrow (S\otimes_k S)_\rho =
{\textstyle\bigoplus_{a=0}^\rho} S_{\rho-a}\otimes_k S_a.
\] 
\end{enumerate}
\end{remark}

Bezoutians can be defined in greater generality (see
\cite{becker,kunz}), but the case considered in
Definition~\ref{bezoutdef} is the only one we need for our main
result.

By Remark~\ref{bezoutprops}, the Bezoutian in degree $(\rho-a,a)$
gives a linear map
\[
{\textstyle\bigwedge^{n+1}} W = {\textstyle\bigwedge^{n+1}} S_d
\longrightarrow S_{\rho-a}\otimes_k S_a,
\] 
which by \eqref{hom} corresponds to an $E$-module homomorphism
\begin{equation}
\label{bpdef}
B_p : \ed(n-p)\otimes_k S_{\rho-a}^* \longrightarrow
\ed(-p-1)\otimes_k S_a.
\end{equation}

\begin{theorem}
\label{main}
The sheaf $\Fc = \nu_{d*}(\Oc_{\P^n}(\ell))$ has a Tate resolution
with 
\[
T^p(\Fc) = \widehat{E}(-p)\otimes_k S_{a-d} \bigoplus
\widehat{E}(n-p)\otimes_k S_{\rho - a}^*, \qquad a = \ell + (p+1)d,
\]
and the differential $d_p : T^p(\Fc) \to T^{p+1}(\Fc)$ is given by
\[
\xymatrix@R=5pt{
\widehat{E}(n-p)\otimes_k S_{\rho - a}^* \ar[r]^(.45){\alpha_p}
\ar[ddr]^{(-1)^p B_p} & \widehat{E}(n-p-1)\otimes_k S_{\rho - a - d}^*\\
\bigoplus & \bigoplus\\
\widehat{E}(-p)\otimes_k S_{a-d} \ar[r]^{\beta_p} &
\widehat{E}(-p-1)\otimes_k S_{a},}
\]
where $B_p$ is the Bezoutian map from \eqref{bpdef} and
$\alpha_p,\beta_p$ are as in \eqref{tpf}.
\end{theorem}

\section{Proof of the Main Result}

We begin with two lemmas needed for the proof of Theorem~\ref{main}.
The notation will be the same as for the previous section.  First
observe that the graded pieces of $B_p$ from \eqref{bpdef} induce
linear maps
\[
{\textstyle\bigwedge^{n+1+m}} W \otimes_k S_{\rho - a}^*
\longrightarrow {\textstyle\bigwedge^{m}} W \otimes_k S_a
\]
for any integer $m$.  This follows from $\ed(n-p)_{p+1+m} =
\bigwedge^{n+1+m} W$.  These maps will be called $B_p$ by abuse of
notation.  Then one of the graded pieces of the differentials $d_p$
from Theorem~\ref{main} give the diagram
\[
\xymatrix@R=5pt{
{\textstyle\bigwedge^{n+2}} W \otimes_k S_{\rho - a}^* \ar[r]^(.48){\alpha_p}
\ar[ddr]^{(-1)^p B_p} &  {\textstyle\bigwedge^{n+1}}W \otimes_k S_{\rho - 
  a - d}^* \ar[ddr]^{(-1)^{p+1}B_{p+1}}  &\\ 
& \bigoplus & \\
& W\otimes_k S_{a} \ar[r]^{\beta_{p+1}} & S_{a+d}.}
\] 

\begin{lemma}
\label{dsquared}
$(-1)^{p+1}B_{p+1}\circ \alpha_p + \beta_{p+1}\circ (-1)^p B_p = 0$ in
the above diagram.
\end{lemma}

\begin{proof} Given $f_0,\dots,f_{n+1} \in W = S_d$, the
polynomials $\Delta_{j}(f_i)$ from Definition~\ref{bezoutdef} satisfy
the identity
\[
\sum_{j=0}^n \Delta_{j}(f_i) (x_i - y_i) = f_i(x) - f_i(y),\quad 0 \le i
\le n+1,
\]
by a telescoping sum argument.  Here we write $f_i(x)$ for
$f_i(x_0,\dots,x_n)$, and similarly for $f_i(y)$.  It follows that in
the $(n+2)\times(n+2)$ matrix
\[
\left(\begin{array}{cccc} 
f_0(x) - f_0(y) & f_1(x) - f_1(y) & \cdots & f_{n+1}(x) - f_{n+1}(y)\\ 
\Delta_{0}(f_0) & \Delta_{0}(f_1) & \cdots & \Delta_{0}(f_{n+1})\\
\vdots          & \vdots          & \ddots & \vdots\\
\Delta_{n}(f_0) & \Delta_{n}(f_1) & \cdots & \Delta_{n}(f_{n+1})
\end{array}\right),
\]
the first row is a linear combination (in $k[x,y]$) of the remaining
rows.  Hence the determinant is zero.  Now expand by minors along the
first row and observe that the $(n+1)\times(n+1)$ minors of the last
$n+1$ rows are Bezoutians.  Hence we get an identity
\[
\sum_{i=0}^{n+1}(-1)^i \Delta^i(x,y) f_i(x) =
\sum_{i=0}^{n+1}(-1)^i \Delta^i(x,y) f_i(y),
\]
where $\Delta^i(x,y)$ is the Bezoutian of
$f_0,\dots,\widehat{f_i},\dots,f_{n+1}$.  Each side is homogeneous of
degree $\rho + d$ in $k[x,y]$, where $\rho = (n+1)(d-1)$.

If we write $\Delta^i(x,y) = \sum_{|\alpha| \le \rho}
\Delta_\alpha^i(x) y^\alpha$, then we can write the identity as
\[
\sum_{i=0}^{n+1}(-1)^i \sum_{|\alpha| \le \rho} \Delta_\alpha^i(x)
f_i(x) y^\alpha =
\sum_{i=0}^{n+1}(-1)^i \sum_{|\alpha| \le \rho} \Delta_\alpha^i(x)
f_i(y) y^\alpha. 
\]
Using $k[x,y] \simeq S\otimes_k S$ and taking the graded piece of
bidegree $(\rho - a,a+d)$ gives
\begin{equation}
\label{identity}
\sum_{i=0}^{n+1}(-1)^i  \sum_{|\alpha| = a+d} \Delta_\alpha^i(x)
f_i(x)\otimes x^\alpha =
\sum_{i=0}^{n+1}(-1)^i \sum_{|\alpha| = a} \Delta_\alpha^i(x)\otimes
f_i(x)x^\alpha. 
\end{equation}
This is an identity in $S_{\rho - a} \otimes_k S_{a+d}$.
 
Now pick $\varphi \in S_{\rho -a}^*$.  If we apply $\varphi \otimes 1$
to \eqref{identity}, we obtain the identity
\begin{equation}
\label{identity2}
\sum_{i=0}^{n+1}  (-1)^i\!\! \sum_{|\alpha| = a+d}\!\! 
\varphi\big(\Delta_\alpha^i(x) f_i(x)\big)\, x^\alpha =
\sum_{i=0}^{n+1}(-1)^i \sum_{|\alpha| = a}
\varphi\big(\Delta_\alpha^i(x)\big)\, f_i(x)x^\alpha
\end{equation} 
in $S_{a+d}$.  The left-hand side of \eqref{identity2} is
$B_{p+1}\circ\alpha_p$ evaluated at $f_0\wedge\cdots\wedge f_{n+1}
\otimes \varphi$, while the right-hand side is $\beta_{p+1} \circ
B_{p}$ evaluated at the same element.  This shows that
$B_{p+1}\circ\alpha_p - \beta_{p+1} \circ B_{p} = 0$, from which the
lemma follows immediately.
\end{proof}

To prepare for the second lemma, let $N = \dim(W) = \binom{n+d}{d}$
and assume that $0 \le \rho - a < d$, so that $S_{\rho-a-d}^* = 0$.
Then one of the graded pieces of the differential $d_p$ from
Theorem~\ref{main} gives the diagram
\begin{equation}
\label{lemma2prep}
\xymatrix@R=5pt{
{\textstyle\bigwedge^{N}} W \otimes_k S_{\rho - a}^* 
\ar[ddr]^{(-1)^p B_p} &   \\ 
\bigoplus & \\
{\textstyle\bigwedge^{N-n}}W\otimes_k S_{a-d} \ar[r]^{\beta_p} &
{\textstyle\bigwedge^{N-n-1}}W\otimes_k S_{a}.}  
\end{equation}

\begin{lemma} 
\label{secondlem}
If $0 \le \rho - a < d$, then the maps $B_p$ and
$\beta_p$ in \eqref{lemma2prep} have the following two properties:
\begin{enumerate}
\item $B_p$ is injective.
\item $\mathrm{Im}(B_p) \cap \mathrm{Im}(\beta_p) = \{0\}$.
\end{enumerate}
\end{lemma}

\begin{proof}
The Bezoutian of $x_0^d,\dots,x_n^d$ is easily seen to be
\[
\Delta = \sum_{\beta \le \beta_{d-1}} x^\beta
y^{\beta_{d-1} -\beta},
\]
where $\beta_{d-1} = (d-1,\dots,d-1) \in \Z^n$ and $\beta \le
\beta_{d-1}$ means that every component of $\beta$ is $\le d-1$.  This
Bezoutian is also computed in \cite{becker}.

The monomial basis of $W = S_d$ induces a basis of $\bigwedge^i W$ for
every $i$.  When $i = N$, the space has dimension one, and we write
its basis element as
\[
x_0^d \wedge \cdots \wedge x_n^d \wedge \omega \in
{\textstyle\bigwedge^N} W,
\]
where $\omega$ is the wedge product of the remaining monomials of
degree $d$.  Given $\varphi \in S_{\rho-a}^*$, we obtain 
\begin{equation}
\label{bpfor}
B_p(x_0^d \wedge \cdots \wedge x_n^d \wedge \omega \otimes \varphi) =
\omega \otimes \big({\textstyle \sum_{\beta}} \varphi(x^\beta)\,
x^{\beta_{d-1}-\beta}\big) + \cdots,
\end{equation}
where the sum inside the parentheses is over all $\beta$ of degree
$\rho - a$ satisfying $\beta \le \beta_{d-1}$, and the omitted terms
involve basis elements of $\bigwedge^{N-n-1}W$ different from
$\omega$.

Let $\varphi$ be in the kernel of $B_p$.  It follows that
$\varphi(x^\beta) = 0$ for all $x^\beta$ appearing in the above sum.
But our hypothesis that $\rho - a < d$ guarantees that this sum
includes \emph{all} monomials of degree $\rho - a$.  These monomials
form a basis of $S_{\rho-a}$, so that $\varphi$ must vanish.  This
proves that $B_p$ is injective, as claimed.

For the second part of the lemma, let $A = \sum_i \omega_i \otimes p_i
\in \bigwedge^{N-n}W\otimes_k S_{a-d}$, where $\{\omega_i\}_i$ is the
basis of $\bigwedge^{N-n}W$ coming from monomials.  We can assume that
the basis includes $\omega_i = \omega \wedge x_i^d $ for $i =
0,\dots,n$, where $\omega$ is as above.  Then
\[
\beta_p(A) = \omega \otimes \big({\textstyle \sum_{i=0}^n} x_i^d
p_i\big) + \cdots,
\]
where the omitted terms involve basis elements of $\bigwedge^{N-n-1}W$
different from $\omega$.  The monomials appearing in $\sum_{i=0}^n
x_i^d p_i$ all have some $x_i$ with an exponent $\ge d$, yet in the
$\omega$-term of \eqref{bpfor}, every $x_i$ has exponent $\le d-1$.
Hence, if $\beta_p(A) = B_p(x_0^d\wedge\cdots\wedge x_n^d\wedge \omega
\otimes \varphi)$, then their $\omega$-terms in $\bigwedge^{N-n-1} W
\otimes_k S_a$ must vanish, which as above implies that $\varphi = 0$.
Hence $\mathrm{Im}(B_p) \cap \mathrm{Im}(\beta_p) = \{0\}$.
\end{proof}

We can now prove our main result.

\begin{proof}[Proof of Theorem~\ref{main}] 
We first show that the differential $d_p : T^p(\Fc) \to T^{p+1}(\Fc)$
defined in Theorem~\ref{main} satisfies $d_{p+1}\circ d_p = 0$, i.e.,
$(T^\bullet(\Fc),d_\bullet)$ is a complex.

We know that $\alpha_{p+1}\circ\alpha_p = 0$ and
$\beta_{p+1}\circ\beta_p = 0$.  It remains to show that the map
\[
\ed(n-p)\otimes S_{\rho-a}^* \longrightarrow \ed(-p-2)\otimes S_{a+d}
\]
given by $(-1)^{p+1}B_{p+1}\circ \alpha_p + \beta_{p+1}\circ (-1)^p
B_p$ is zero.  Since
\[
\mathrm{Hom}_E(\ed(n-p)\otimes S_{\rho-a}^*,\ed(-p-2)\otimes
S_{a+d})_0 \simeq \mathrm{Hom}_k({\textstyle\bigwedge^{n+2}} W\otimes
S_{\rho-a}^*,S_{a+d}),
\]
this follows immediately from Lemma~\ref{dsquared}.

Next we need to show that for each $p$, $d_p$ is determined by the
minimal generators of the kernel of $d_{p+1}$.  This is where we use
the power of the formula for $T^p(\Fc)$ given in \eqref{tp}:\ it tells
us the degrees of the minimal generators of $\mathrm{Ker}(d_{p+1})$
and the number of minimal generators in these degrees.  Furthermores,
$d_{p+1}\circ d_p = 0$ implies that $d_p$ maps into the kernel.  So we
need to study how $d_p$ behaves in the degrees of the minimal
generators.

Recall that $a = \ell + (p+1)d$, so that $\rho - a < 0$ for large $p$.
We will look closely at the case when $0 \le \rho - a < d$.  Here,
$d_{p+1} = \beta_{p+1}$ and the complex looks like
\[
\xymatrix@R=5pt{
\widehat{E}(n-p)\otimes_k S_{\rho - a}^* 
\ar[ddr]^{(-1)^p B_p} & &\\
\bigoplus & & \\
\widehat{E}(-p)\otimes_k S_{a-d} \ar[r]^{\beta_p} &
\widehat{E}(-p-1)\otimes_k S_{a} \ar[r]^(.47){\beta_{p+1}} &
\widehat{E}(-p-2)\otimes_k S_{a+d}.} 
\]
This is the first place where a nonzero diagonal map appears in the
Tate resolution.  Since $\ed \simeq E(-N)$ (this is the notation of
Lemma~\ref{secondlem}), there are $\dim(S_{a-d})$ minimal generators
of degree $N+p$ and $\dim(S_{\rho - a}^*)$ minimal generators of
degree $N-n+p$.  The former are taken care of by the known formula for
$\beta_p$.  For the latter, notice that the above diagram in degree
$N-n+p$ is precisely \eqref{lemma2prep}, and then
Lemma~\ref{secondlem} implies that $(-1)^p B_p$ maps injectively onto
the minimal generators in this degree.  Hence we have the desired
behavior when $\rho - a < d$.

We now proceed by decreasing induction on $p$.  Suppose that $\rho - a
\ge d$ and that everything is fine for larger $p$.  As above, there
are $\dim(S_{a-d})$ minimal generators of degree $N+p$ and
$\dim(S_{\rho - a}^*)$ minimal generators of degree $N-n+p$, where the
former are taken care of by $\beta_p$.  But now in degree $N-n+p$,
the differential $d_p$ is given by
\[
\xymatrix@R=5pt{
{\textstyle\bigwedge^N} W \otimes_k S_{\rho - a}^* \ar[r]^(.45){\alpha_p}
\ar[ddr]^{(-1)^p B_p} & {\textstyle\bigwedge^{N-1}} W\otimes_k 
S_{\rho - a - d}^*\\ \bigoplus & \bigoplus\\
{\textstyle\bigwedge^{N-n}} W\otimes_k S_{a-d} \ar[r]^{\beta_p} &
{\textstyle\bigwedge^{N-n-1}} W\otimes_k S_{a}.}
\]
The key observation is that the $\alpha_p$ in this diagram is dual to
the multiplication map $W\otimes S_{\rho-a-d} \to S_{\rho-a}$, which
is surjective since $\rho-a \ge d$.  This implies that in the degree
of the minimal generators, $\alpha_p$ is injective.  It follows that
$\alpha_p\oplus(-1)^p B_p$ is injective in this degree and its image
intersects the image of $\beta_p$ in $\{0\}$.  This shows that $d_p$
has the desired property and completes the proof of the theorem.
\end{proof}

\begin{remark} 
As noted by Materov, the Tate resolution of Theorem~\ref{main} can be
expressed as a mapping cone.  Let $\mathcal{C}^\bullet$ be the part of
the Tate resolution in cohomological degree $n$ (i.e., the part of
\eqref{tp} involving $H^n$).  Thus $\mathcal{C}^\bullet$ is given by
\[
\cdots \longrightarrow \mathcal{C}^p = \ed(n-p) \otimes_k S_{a-d}^*
\stackrel{\alpha_p}{\longrightarrow} \mathcal{C}^{p+1} =
\ed(n-p-1)\otimes_k S_{a} \longrightarrow \cdots.
\]
Similarly, let $\mathcal{D}^\bullet$ denote the part of the Tate
resolution in cohomological degree 0, shifted by 1.  Thus
$\mathcal{D}^\bullet$ is given by
\[
\cdots \longrightarrow \mathcal{D}^p = \ed(-p-1) \otimes_k S_{a}
\xrightarrow{\beta_{p+1}} \mathcal{D}^{p+1} =
\ed(-p-2)\otimes_k S_{a+d} \longrightarrow \cdots.
\]
The proofs of Lemma~\ref{dsquared} and Theorem~\ref{main} give a
commutative diagram
\[
\xymatrix{
\cdots \ar[r] & \ed(n-p)\otimes_k S^*_{\rho - a} \ar[r]^(.45){\alpha_p}
\ar[d]^{B_p} & \ed(n-p-1)\otimes_k S^*_{\rho - a-d} \ar[r]
\ar[d]^{B_{p+1}} &\cdots \\
\cdots \ar[r] & \ed(-p-1)\otimes_k S_{a} \ar[r]^(.47){\beta_{p+1}} &
\ed(-p-2)\otimes_k S_{a+d} \ar[r] & \cdots,
}
\]
so that the Bezoutians $\{B_p\}$ give a map of complexes
$\mathcal{C}^\bullet \to \mathcal{D}^\bullet$.  Then
Theorem~\ref{main} implies that the Tate resolution is the mapping
cone of this map of complexes.  This explains the signs $(-1)^p$ and
$(-1)^{p+1}$ appearing in the statement of the theorem.
\end{remark}

\section{Application to Duality}

We conclude by exploring the relation between duality, Bezoutians, and
the Tate resolution.  We first recall how to extract information from
the Tate resolution.  Stated briefly, the key idea is to look at
$T^\bullet(\Fc)$ in a specific degree, but only \emph{after} replacing
$W$ with a suitable subspace $U \subset W$.  This is the functor
$\mathbf{U}_l$ from \cite{es}, which is equivalent to the projection
formula from \cite[Sect.\ 1.2]{floystad}.

To make this precise, let $U \subset W$ be a subspace.  Since $\P(W) =
(W^*-\{0\})/k^*$, the linear subspace $\P(W/U) \subset \P(W)$ is the
center of the projection $\pi : \P(W) \dashrightarrow \P(U)$.  If
$\P(W/U)$ is disjoint from the support of $\Fc$, then \cite{es} and
\cite{floystad} show that
\[
T_U^\bullet(\Fc) = \mathrm{Hom}_E({\textstyle\bigwedge}
U^*,T^\bullet(\Fc))
\]
is a Tate resolution of $\pi_*\Fc$ on $\P(U)$.  Note also that $\Fc$
and $\pi_*\Fc$ have the same cohomology since $\pi : \P(W) \setminus
\P(W/U) \rightarrow \P(U)$ is affine.

In the situation of Theorem~\ref{main}, we have $W = S_d$, so that a
subspace $U \subset W$ satisfies
\[
\P(W/U) \cap \mathrm{Supp}(\Fc) = \emptyset
\]
if and only if the homogeneous polynomials in $U$ have no common zeros
in $\P^n$.  When this happens, the above paragraph and
Theorem~\ref{main} give a minimal exact sequence of free graded
$E_U$-modules $T_U^\bullet(\Fc)$, where $T_U^p(\Fc) \to
T_U^{p+1}(\Fc)$ is
\begin{equation}
\label{tru}
\xymatrix@R=5pt{
\widehat{E}_U(n-p)\otimes_k S_{\rho - a}^* \ar[r]^(.45){\alpha_p}
\ar[ddr]^{(-1)^p B_p} & \widehat{E}_U(n-p-1)\otimes_k S_{\rho - a - d}^*\\
\bigoplus & \bigoplus\\
\widehat{E}_U(-p)\otimes_k S_{a-d} \ar[r]^{\beta_p} &
\widehat{E}_U(-p-1)\otimes_k S_{a}.}
\end{equation}
Here, $E_U = \bigwedge U^*$ and $\ed_U = \bigwedge U$.  As we will
see, looking at this complex in specific degrees for specific choices
of $U$ will give some interesting duality theorems.

\begin{example}
\label{first}
First let $U = \mathrm{Span}(f_0,\dots,f_{n}) \subset W = S_d$,
where $f_0,\dots,f_n$ have no common zeros on $\P^n$.  As is
well-known, this happens $\iff$ $f_0,\dots,f_n$ is a regular sequence
$\iff$ the Koszul complex of $f_0,\dots,f_n$ is exact.

Let $I = \langle f_0,\dots,f_n\rangle \subset S$ and $R = S/I$.  Then
consider $T_U^\bullet(\Fc)$ in degree $p+1$.  Using \eqref{tru}, we
obtain the following exact sequence of vector spaces:
\[
\xymatrix@R=5pt{
& {\textstyle\bigwedge^{n+1}} U \otimes_k S_{\rho - a}^* \ar[r]^{\alpha_p}
\ar[ddr]^{(-1)^p B_p} & {\textstyle\bigwedge^{n}} U\otimes_k 
S_{\rho - a - d}^* \ar[r] & \cdots\\ 
&\bigoplus & \bigoplus &\\
\cdots \ar[r] & U\otimes_k S_{a-d} \ar[r]^{\beta_p} & S_{a}. & }
\]
It follows that $(-1)^p B_p$ induces an isomorphism
\[
\mathrm{Ker}(\alpha_p) \simeq \mathrm{Coker}(\beta_p).
\]
Since $\mathrm{Ker}(\alpha_p) = R_{\rho-a}^*$ and
$\mathrm{Coker}(\beta_p) = R_a$, we recover the known duality
\[
R_{\rho -a}^* \simeq R_a.
\]
Furthermore, $\bigwedge^{n+1} U$ has basis element
$f_0\wedge\cdots\wedge f_n$, so that if 
\[
\Delta = \sum_{|\alpha| \le \rho} \Delta_\alpha(x)y^\alpha
\]
is the Bezoutian of $f_0,\dots,f_n$, then the above isomorphism
$R_{\rho -a}^* \simeq R_a$ is given by
\begin{equation}
\label{bezduality1}
\varphi \in R_{\rho-a}^* \longmapsto \sum_{|\alpha| = a}
\varphi([\Delta_\alpha(x)])\, [x^\alpha] \in R_a,
\end{equation}
where $[g] \in R$ denotes the coset of the polynomial $g \in S$.
\end{example}

\begin{remark} Here are some comments about Example~\ref{first}.
\begin{enumerate}
\item It is known that the duality $R_{\rho -a}^* \simeq R_a$ can be
computed by \eqref{bezduality1}.  Proofs can be found in
\cite{becker,kunz} in the case when the $f_i$ are homogeneous of
degree $d_i$, as opposed to the equal degree case considered here.
Our contribution is to show that the Tate resolution gives a new proof
of this explicit duality in the equal degree case.
\item The proof given in \cite{becker} that \eqref{bezduality1}
induces $R_{\rho -a}^* \simeq R_a$ uses the Bezoutian of
$x_0^d,\dots,x_n^d$.  This is the same Bezoutian used in the proof of
Lemma~\ref{secondlem}.
\end{enumerate}
\end{remark}

\begin{example}
\label{second}
Now suppose that $U = \mathrm{Span}(f_0,\dots,f_n,f_{n+1}) \subset W$,
where the polynomials $f_0,\dots,f_n,f_{n+1}$ are linearly independent
and have no common zeros in $\P^n$.  We have one more polynomial than
we had in Example~\ref{first}.  As we will see, this leads to a
slightly different form of duality.

As in the previous example, let $I = \langle
f_0,\dots,f_n,f_{n+1}\rangle \subset S$ and $R = S/I$, and consider
$T_U^\bullet(\Fc)$ in degree $p+2$.  Using \eqref{tru}, we obtain the
following exact sequence of vector spaces:
\[
\xymatrix@R=5pt@C=15pt{ 
& {\textstyle\bigwedge^{n+2}} U \otimes_k S_{\rho - a}^*
  \ar[r]^(.47){\alpha_p} \ar[ddr]^{(-1)^p B_p} & 
{\textstyle\bigwedge^{n+1}} U\otimes_k S_{\rho - a - d}^*
\ar[r]^{\alpha_{p+1}} \ar[ddr]^{(-1)^{p+1} B_{p+1}}
& {\textstyle\bigwedge^{n}} U\otimes_k S_{\rho - a - 2d}^* \ar[r] &\cdots\\ 
&\bigoplus & \bigoplus & \bigoplus & \\ \cdots \ar[r] &
{\textstyle\bigwedge^{2}}U\otimes_k S_{a-d} \ar[r]^{\beta_p} &
U \otimes_k S_{a} \ar[r]^{\beta_{p+1}} & S_{a+d}.}
\]
It follows that $(-1)^p B_p$ induces an isomorphism
\begin{equation}
\label{duality21}
\mathrm{Ker}(\alpha_p) \simeq
\mathrm{Ker}(\beta_{p+1})/\mathrm{Im}(\beta_p).
\end{equation}
Note that $\mathrm{Ker}(\alpha_p) = R_{\rho-a}^*$ and that the bottom
row of the above diagram comes from the Koszul complex of
$f_0,\dots,f_{n+1}$.  Hence
\[
\mathrm{Ker}(\beta_{p+1}) = \mathrm{Syz}(f_0,\dots,f_{n+1})_{a+d},
\]
where a syzygy $(A_0,\dots,A_{n+1})$ is said to have degree $a+d$ if
$\sum_{i=0}^{n+1} A_i f_i = 0$ in $S_{a+d}$.  Furthermore, the image
of $\beta_p : \bigwedge^{2} U\otimes_k S_{a-d} \to U \otimes_k S_a$
is the submodule of $\mathrm{Syz}(f_0,\dots,f_{n+1})_{a+d}$ consisting
of Koszul syzygies.  Hence we set
\[
\mathrm{Kosz}_{a+d} = \mathrm{Im}(\beta_p).
\]
Then the duality \eqref{duality21} becomes
\begin{equation}
\label{duality22}
R_{\rho - a}^* \simeq
\mathrm{Syz}(f_0,\dots,f_{n+1})_{a+d}/\mathrm{Kosz}_{a+d}.
\end{equation}
Notice also that $B_p$ gives an explicit description of this duality
since elements of $R_{\rho - a}^*$ can be regarded as linear
functionals $\varphi$ on $S_{\rho - a}$ that vanish on $I_{\rho-a}$.
Then the left-hand side of \eqref{identity2} vanishes, so that
\eqref{identity2} becomes
\begin{equation}
\label{koszulsyz}
\sum_{i=0}^{n+1} (-1)^i \sum_{|\alpha| = a}
\varphi\big(\Delta^i_{\alpha}\big)x^\alpha f_i = 0.
\end{equation}
As noted in the proof of Lemma~\ref{secondlem}, this is
$\beta_{p+1}$ applied to $B_p(f_0\wedge\cdots\wedge f_{n+1}\otimes
\varphi)$.  Thus
\[
\Big({\textstyle \sum_{|\alpha| = a}}
\varphi\big(\Delta^0_{\alpha}\big)x^\alpha, 
-{\textstyle \sum_{|\alpha| = a}}
\varphi\big(\Delta^1_{\alpha}\big)x^\alpha,\dots,
(-1)^{n+1}{\textstyle \sum_{|\alpha| = a}}
\varphi\big(\Delta^{n+1}_{\alpha}\big)x^\alpha\Big)
\]
is an element of $\mathrm{Syz}(f_0,\dots,f_{n+1})_{a+d}$ coming from
$B_p$.  We call this a \emph{Bezout syzygy}.  It follows that the
duality \eqref{duality22} is computed in terms of Bezout syzygies.
\end{example}

\begin{remark} 
\label{final}
Here are further comments on the duality of
Example~\ref{second}. 
\begin{enumerate}
\item If $K_\bullet$ is the Koszul complex of $f_0,\dots,f_{n+1}$,
  then our hypothesis that the $f_i$ don't vanish simultaneously on
  $\P^n$ implies that $K_\bullet$ is almost exact.  In fact, the only
  place exactness fails is at $K_1$:
\[
\begin{array}{ccccccccccc}
\cdots & \longrightarrow & K_2 & \stackrel{d_1}{\longrightarrow} & K_1
& \stackrel{d_0}{\longrightarrow} & 
S & \longrightarrow & R & \longrightarrow & 0\\
\uparrow && \uparrow && \uparrow && \uparrow && \uparrow &&\\
\mbox{\small ok} && \mbox{\small ok} && \mbox{\small \bf no} &&
\mbox{\small ok} && \mbox{\small ok} && 
\end{array}
\] 
The graded pieces of $\mathrm{Ker}(d_0)/\mathrm{Im}(d_1)$ are the
$\mathrm{Syz}(f_0,\dots,f_{n+1})_{a+d}/\mathrm{Kosz}_{a+d}$ appearing
in \eqref{duality22}.  Thus size of
\[
R = S/I = k[x_0,\dots,x_{n}]/\langle f_0,\dots,f_{n+1}\rangle
\]
gives a precise measure of the failure of an arbitrary syzygy to be
Koszul. 
\item One corollary of the duality \eqref{duality22} is that the
syzygy module of $f_0,\dots,f_{n+1}$ is generated by Koszul syzygies
and Bezout syzygies. 
\item We can write the duality \eqref{duality22} more conceptually as
follows.  Set $\sigma = \sum_{i=0}^{n+1} \deg(f_i) - (n+1) = \rho + d$
and $b = a+d$.  Then \eqref{duality22} becomes 
\[
R_{\sigma - b}^* \simeq
\mathrm{Syz}(f_0,\dots,f_{n+1})_{b}/\mathrm{Kosz}_{b}.
\]
Furthermore, if $H_i(K_\bullet)$ is the $i$th homology of the Koszul
complex, then this duality can be written as
\[
H_0(K_\bullet)_{\sigma - b}^* \simeq H_1(K_\bullet)_b.
\]
\item More generally, suppose that $f_0,\dots,f_m \in S_d$ are
linearly independent and don't vanish simultaneously on $\P^n$.  Note
that $m \ge n$ and that Examples~\ref{first} and \ref{second}
correspond to $m = n$ and $m = n+1$ respectively.  Let $K_\bullet$ be
the Koszul complex of $f_0,\dots,f_m$ and set $\sigma = \sum_{i=0}^{m}
\deg(f_i) - (n+1)$.  Then Examples~\ref{first} and \ref{second} easily
generalize to give a \emph{Kozsul duality}
\[
H_i(K_\bullet)^*_{\sigma-a} \simeq H_{m-n-i}(K_\bullet)_a,\qquad 0 \le
i \le m-n,
\]
that is computed by Bezoutians.
\item The Koszul duality just stated applies more generally to
homogeneous polynomials in $S$ of arbitrary degrees (not necessarily
equal) that don't vanish simultaneously on $\P^n$.  The proof that
some isomorphism exists is an easy spectral sequence argument; the
fact that it is given by Bezoutians takes more work---this has been
proved by Jouanolou (unpublished).  So again, the Tate resolution
gives a quick proof of the equal degree case of an explicit duality
theorem.
\end{enumerate}
\end{remark}

A final comment is that the duality theorems of Examples~\ref{first}
and \ref{second} and Remark~\ref{final} come from the \emph{same} Tate
resolution.  Once we describe the Tate resolution in terms of
Bezoutians, we get immediate Bezoutian descriptions of \emph{all} of
these duality results.  This indicates the deep relation between
duality, Bezoutians, and the Tate resolution.

\section*{Acknowledgements}

I would like to thank Evgeny Materov, Jessica Sidman and Ivan
Soprounov for their useful observations and careful reading of the
paper.  Thanks also go to Laurent Bus\'e, Marc Chardin, David Eisenbud
and Jean-Pierre Jouanolou for helpful conversations.

\end{document}